\documentclass[a4paper]{article}
\usepackage[T1]{fontenc}
\usepackage{latexsym}
\usepackage{amssymb}
\usepackage{cite}

\newcommand{\Ff}{{\mathbb F}}

\newcommand{\adj}{\,\sim\,}
\newcommand\textprime{$'$}

\title{Some locally Kneser graphs}
\author{A. E. Brouwer}
\date{2023-12-05}

\begin{document}
\maketitle

\begin{abstract}
The Kneser graph $K(n,d)$ is the graph on the $d$-subsets of an $n$-set,
adjacent when disjoint. Clearly, $K(n+d,d)$ is locally $K(n,d)$.
Hall showed for $n \ge 3d+1$ that there are no further examples.
Here we give other examples of locally $K(n,d)$ graphs for $n = 3d$,
and some further sporadic examples. It follows that Hall's bound is
best possible.
\end{abstract}

\section{Locally something graphs}
A graph $\Gamma$ is called {\em locally $\Delta$} when for each vertex
of $\Gamma$ the subgraph induced on the set of its neighbors is isomorphic
to $\Delta$.
The Trahtenbrot-Zykov problem~\cite{Z63} asks whether given a finite graph
$\Delta$ there exists a graph $\Gamma$ that is locally $\Delta$.
In general, this question is undecidable (Bulitko \cite{B73}).
It is unknown whether the problem restricted to finite $\Gamma$
is also undecidable.

More generally, one wants to classify all such graphs $\Gamma$.
Hall \cite{H85} determines the possible $\Gamma$ for all graphs $\Delta$
on at most 6 vertices.
For some $\Delta$
a graph $\Gamma$ that is locally $\Delta$ is necessarily infinite.

Weetman \cite{W94,W94a} constructs infinite locally $\Delta$ graphs
for $\Delta$ of girth at least 6, and proves a diameter bound
in certain other cases.

There is a large literature, see e.g. \cite{BHM80,BB89,BB92,BBBC84,
BC73,BC75,BH77,BCN,BKK03,BvM,Bug90,Bug92,B72,B73,B83,H80,H85,H87,
HS85,H78,LM13,N93,R81,vdC85,V81,W94,W94a,Y21}.

\section{Locally Kneser graphs}
The {\em Kneser graph} $K(n,d)$ (where $0 \le d \le n$) is the graph on
the $d$-subsets of an $n$-set, adjacent when disjoint.

\medskip
Hall \cite{H87} shows that for $n \ge 3d+1$ any connected
locally $K(n,d)$ graph is isomorpic to $K(n+d,d)$, and wonders
whether this bound can be improved.

In fact there are further examples for $n = 3d$:
The graph on the $2^{n-1}$ even weight binary vectors of length $n$,
adjacent when their difference has weight $2d$ is locally
$K(n,d)$ and different from $K(n+d,d)$ (for $d > 1$).
It follows that Hall's bound is best possible.

\medskip
For $n = 2d+1$, $d \ge 3$, the graph $K(n,d)$ has girth 6,
so that there exist infinite locally $K(n,d)$ graphs
by Weetman \cite{W94}.

\medskip
There are three locally $K(5,2)$ graphs, on 21, 63, and 65 vertices
(Hall \cite{H80}).
Note that $K(5,2)$ is the Petersen graph.

\medskip
There are three locally $K(6,2)$ graphs, on 28, 32, and 36 vertices
(\-Buekenhout \& Hubaut \cite{BH77}).
Note that $K(6,2)$ is the collinearity graph of the generalized quadrangle
of order 2.

\medskip
As we saw, there are infinite locally $K(7,3)$ graphs.
(Are there also finite examples?)

\medskip
The graph on the elliptic lines in the $O_8^-(2)$ geometry,
adjacent when orthogonal, is locally $K(8,3)$.
The automorphism group is $O_8^-(2){:}2$, with point stabilizer
$S_3 \times S_8$. Diagram:

{\scriptsize

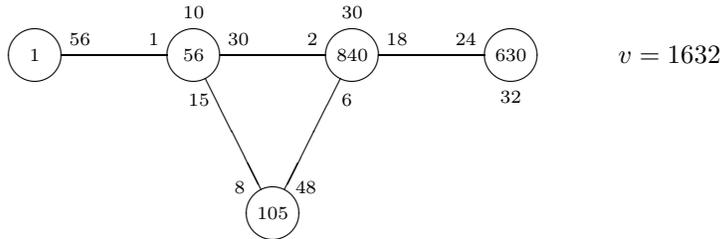
$$\begin{picture}(240,80)(0,-60)
\put(0,0){\circle{20}}\put(0,0){\makebox(0,0){1}}
\put(10,0){\line(1,0){40}}
\put(13,6){\makebox(0,0)[l]{56}}
\put(47,6){\makebox(0,0)[r]{1}}
\put(60,0){\circle{20}}\put(60,0){\makebox(0,0){56}}
\put(60,16){\makebox(0,0){10}}
\put(70,0){\line(1,0){40}}
\put(73,6){\makebox(0,0)[l]{30}}
\put(107,6){\makebox(0,0)[r]{2}}
\put(120,0){\circle{20}}\put(120,0){\makebox(0,0){840}}
\put(120,16){\makebox(0,0){30}}
\put(130,0){\line(1,0){40}}
\put(133,6){\makebox(0,0)[l]{18}}
\put(167,6){\makebox(0,0)[r]{24}}
\put(180,0){\circle{20}}\put(180,0){\makebox(0,0){630}}
\put(180,-16){\makebox(0,0){32}}
\put(90,-60){\circle{20}}\put(90,-60){\makebox(0,0){105}}
\put(64.47,-8.94){\line(1,-2){21.06}}
\put(115.53,-8.94){\line(-1,-2){21.06}}
\put(62,-17){\makebox(0,0){15}}
\put(118,-17){\makebox(0,0){6}}
\put(77.5,-50){\makebox(0,0){8}}
\put(102.5,-50){\makebox(0,0){48}}
\put(90,-76){\makebox(0,0){-}}
\put(240,0){\makebox(0,0){\normalsize $v=1632$}}
\end{picture}$$
}

See also \cite{C93}.

\medskip
As we saw, the graph on the 256 binary even weight vectors of length 9,
adjacent when they have distance 6, is locally $K(9,3)$.
Diagram:

{\scriptsize

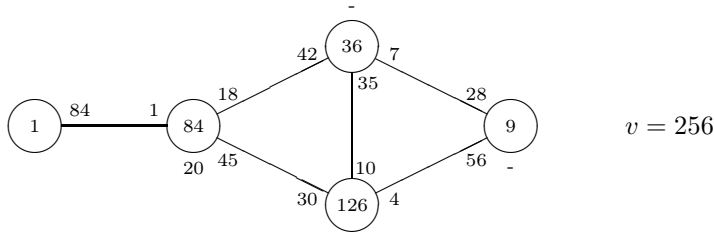
$$\begin{picture}(240,80)(0,-30)
\put(0,0){\circle{20}}\put(0,0){\makebox(0,0){1}}
\put(10,0){\line(1,0){40}}
\put(13,6){\makebox(0,0)[l]{84}}
\put(47,6){\makebox(0,0)[r]{1}}
\put(60,0){\circle{20}}\put(60,0){\makebox(0,0){84}}
\put(60,-16){\makebox(0,0){20}}
\put(68.94430,4.47215){\line(2,1){42.1114}}
\put(68.94430,-4.47215){\line(2,-1){42.1114}}
\put(73,12){\makebox(0,0){18}}
\put(73,-13){\makebox(0,0){45}}
\put(103,27){\makebox(0,0){42}}
\put(103,-28){\makebox(0,0){30}}
\put(120,30){\circle{20}}\put(120,30){\makebox(0,0){36}}
\put(120,45){\makebox(0,0){-}}
\put(120,-46){\makebox(0,0){40}}
\put(126,16){\makebox(0,0){35}}
\put(125,-16){\makebox(0,0){10}}
\put(128.94430,25.52785){\line(2,-1){42.1114}}
\put(128.94430,-25.52785){\line(2,1){42.1114}}
\put(136,27){\makebox(0,0){7}}
\put(136,-28){\makebox(0,0){4}}
\put(167,12){\makebox(0,0){28}}
\put(167,-13){\makebox(0,0){56}}
\put(180,0){\circle{20}}\put(180,0){\makebox(0,0){9}}
\put(180,-16){\makebox(0,0){-}}
\put(120,-30){\circle{20}}\put(120,-30){\makebox(0,0){126}}
\put(120,20){\line(0,-1){40}}
\put(240,0){\makebox(0,0){\normalsize $v=256$}}
\end{picture}$$
}

\medskip
It can be shown using the arguments from \cite{W94a}
that a locally $K(9,3)$ graph is necessarily finite.
(Is the same true for all locally $K(3d,d)$ graphs?)

\medskip
The bimonster $G = M \,{\rm wr}\, 2$ (where $M$ is the monster)
contains a $S_5$-subgroup $S$ whose centralizer is a subgroup $S_{12}$
in which a 7-point stabilizer is conjugate to $S$, see \cite{CP90}.
Let $\Gamma$ be the graph on the $S_5$-subgroups of $G$ conjugate to $S$,
adjacent when they commute. Then $\Gamma$ is locally $K(12,5)$.

\section{Locally $\lambda = 1$ graphs}
Let $\Delta$ be a graph in which every edge is in a unique triangle
(so that $\Delta$ is the collinearity graph of a partial linear space
with lines of size 3). The Kneser graphs $K(3d,d)$ are examples of such graphs.
We study locally $\Delta$ graphs. The special case where $\Delta$
is the line graph of the Petersen graph was studied in \cite{BKK03}.

Given a partial linear space $(X,L)$ with lines of size 3,
let $G$ be the group
$$
G = \langle X \mid x^2 = 1 = xyz ~{\rm for~all}~ x \in X
~{\rm and}~ \{x,y,z\} \in L \rangle .
$$

{\small
Suppose $A,B$ are two disjoint hyperplane complements in $(X,L)$,
so that each line meets $A$ and $B$ in 0 or 2 points.
Then the map that sends the elements of $A$, $B$, $X \setminus (A \cup B)$
to $a$, $b$, and 1, respectively, is a map from $G$ onto the infinite
group $\langle a, b \mid a^2 = b^2 = 1 \rangle$, so that $G$ is infinite.
\par}

\medskip
Given a subgroup $H$ of $G$, let $\Gamma = \Gamma(G,H,X)$
be the graph that has as vertices the cosets $gH$ for $g \in G$,
and adjacencies $g_1H \adj g_2H$ when $g_2^{-1}g_1 \in HXH$.
Assume that $H$ is normal in $G$.
Now the neighbours of $gH$ are the vertices $gxH$ for $x \in X$.
The group $G$ acts vertex transitively on $\Gamma$. The local graph
induced on the set of neighbours of the vertex $H$ of $\Gamma$
has vertex set $\{ xH \mid x \in X \}$, and if $\{x,y,z\}$ is a line,
then $xy = z$ in $G$, so that $xH \adj yH$.
It follows that $x \mapsto xH$ is a homomorphism from the collinearity graph
$\Delta$ of $(X,L)$ onto the $\Gamma$-neighbourhood of $H$.

Is this map injective? Suppose $H$ is contained in the commutator subgroup
$G'$ of $G$. Then $\Gamma(G,H,X)$ has quotient $\Gamma(G,G',X)$ and the
latter can be identified with the Cayley graph with difference set $X$
in the $\Ff_2$-vector space $\langle X \mid x+y+z=0 ~{\rm for}~
\{x,y,z\} \in L \rangle$. If $N$ is the point-line incidence matrix
of $(X,L)$, then this is $\langle X \rangle / N \langle L \rangle$,
where cosets at Hamming distance 1 are adjacent.
Two points remain distinct in the quotient if the column space of $N$
does not contain vectors of weight 2. No additional adjacencies are
introduced if the columns of $N$ are the only vectors of weight 3
in the column space of $N$.

For these latter two conditions to hold, it suffices that for any two
distinct points, and for any three pairwise noncollinear points,
$(X,L)$ has a geometric hyperplane missing precisely one of these points.
%

If these conditions hold, $\Gamma$ is locally $\Delta$.

\subsection*{Examples}
We construct $\Gamma = \Gamma(G,G',X)$ for a number of spaces $(X,L)$
with lines of size~3. Note that $V = G/G'$ is a binary vector space.
In all cases except d), the graph $\Gamma$ is locally $\Delta$,
where $\Delta$ is the collinearity graph of $(X,L)$.

\medskip\noindent
\begin{tabular}{@{}c@{~}c@{~}c@{~~}c@{~}c@{}c@{~~}c@{~~}c@{~~}c@{~~}c@{}}
 & $\Delta$ & parameters & $|G|$ & $\dim\,V$ & $v$ & $k$ & $d$ & rk & ${\rm Aut}\,\Gamma$ \\
\hline
a) & $GQ(2,1)$ & ${\rm srg}(9,4,1,2)$ & 16 & 4 &
  16 & 9 & 2 & 3 & $[2^7 . 3^2]$ \\
b) & $GQ(2,2)$ & ${\rm srg}(15,6,1,3)$ & 32 & 5 &
  32 & 15 & 3 & 4 & $2^5{:}S_6$ \\
c) & $GQ(2,4)$ & ${\rm srg}(27,10,1,5)$ & 64 & 6 &
  64 & 27 & 2 & 3 & $2^6{:}O_6^-(2)$ \\
d) & $VO_4^-(3)$ & ${\rm srg}(81,20,1,6)$ & 1 & 0 & - \\
e) & $L(K(5,2))$ & $\{4,2,1;\,1,1,4\}$ & $\infty$ & 6 &
  64 & 15 & 3 & 6 & $2^6{:}S_5$ \\
f) & $GH(2,1)$ & $\{4,2,2;\,1,1,2\}$ & $\infty$ & 8 &
  256 & 21 & 3 & 7 & $2^8{:}PGL(3,2)$ \\
g) & $GH(2,2)$ & $\{6,4,4;\,1,1,1\}$ & $\infty$ & 14 &
  16384 & 63 & 4 & 15 & $2^{14}{:}G_2(2)$ \\
g\textprime) & $GH(2,2)$ & $\{6,4,4;\,1,1,1\}$ & $\infty$ & 14 &
  16384 & 63 & 6 & 26 & $2^{14}{:}G_2(2)$ \\
h) & $3^3$ & $\{6,4,2;\,1,2,3\}$ & 512 & 8 &
  256 & 27 & 3 & 6 & $[2^{12} . 3^4]$ \\
i) & $3^4$ & $\{8,6,4,2;\,1,2,3,4\}$ & & 16 &
  65536 & 81 & 6 & 30 & $[2^{23} . 3^5]$ \\
j) & $GO(2,1)$ & $\{4,2,2,2;\,1,1,1,2\}$ & $\infty$ & 16 &
  65536 & 45 & 6 & 93 & $2^{16}{:}M_{10}$ \\
k) & $3S_6$ & $\{6,4,2,1;\,1,1,4,6\}$ & $\infty$ & 11 &
  2048 & 45 & 5 & 16 & $2^5{:}(2^6{:}3.S_6)$ \\
l) & $K(9,3)$ & $(v, k)_\Delta=(84, 20)$ & 256 & 8 &
  256 & 84 & 3 & 5 & $2^8{:}S_9$ \\
m) & $K(12,4)$ & $(v, k)_\Delta=(495, 70)$ & 2048 & 11 &
  2048 & 495 & 3 & 7 & $2^{11}{:}S_{12}$
\end{tabular}

\medskip
We give the parameters for a strongly regular graph
as ${\rm srg}(v,k,\lambda,\mu)$, and for a distance-regular graph
of diameter at least 3 as
$\{b_0,\ldots,b_{d-1};\,c_1,\ldots,c_d\}$ (cf.~\cite{BCN}).
In cases a), b), the graphs $\Delta$ are $3^2$ and $K(6,2)$. 
In cases a), b), and c), the graphs $\Gamma$ are $VO_4^+(2)$,
$T\Delta$, and $VO_6^-(2)$ (with notation as in \cite{BvM}).


Cases g) and g\textprime) are the dual Cayley and the Cayley generalized hexagon,
respectively. For the former $\Gamma$ has diameter 4 and
${\rm Aut}\,\Gamma$ has trivial center, for the latter
$\Gamma$ is antipodal of diameter 6 and ${\rm Aut}\,\Gamma$
has a center of order 2 that interchanges antipodes.

In case h) the diagram is

{\scriptsize

$$\begin{picture}(240,80)(0,-30)
\put(0,0){\circle{20}}\put(0,0){\makebox(0,0){1}}
\put(10,0){\line(1,0){40}}
\put(13,6){\makebox(0,0)[l]{27}}
\put(47,6){\makebox(0,0)[r]{1}}
\put(60,0){\circle{20}}\put(60,0){\makebox(0,0){27}}
\put(60,16){\makebox(0,0){6}}
\put(68.94430,4.47215){\line(2,1){42.1114}}
\put(68.94430,-4.47215){\line(2,-1){42.1114}}
\put(73,12){\makebox(0,0){12}}
\put(73,-13){\makebox(0,0){8}}
\put(103,27){\makebox(0,0){6}}
\put(103,-28){\makebox(0,0){2}}
\put(120,30){\circle{20}}\put(120,30){\makebox(0,0){54}}
\put(120,-30){\circle{20}}\put(120,-30){\makebox(0,0){108}}
\put(120,20){\line(0,-1){40}}
\put(120,45){\makebox(0,0){3}}
\put(120,-46){\makebox(0,0){12}}
\put(124,-16){\makebox(0,0){6}}
\put(126,16){\makebox(0,0){12}}
\put(127.071,-22.929){\line(1,1){45.858}}
\put(130,30){\line(1,0){40}}
\put(130,-30){\line(1,0){40}}
\put(133,-23){\makebox(0,0){6}}
\put(135,36){\makebox(0,0){6}}
\put(135,-36){\makebox(0,0){1}}
\put(166,36){\makebox(0,0){6}}
\put(166,-36){\makebox(0,0){9}}
\put(172,16){\makebox(0,0){12}}
\put(180,30){\circle{20}}\put(180,30){\makebox(0,0){54}}
\put(180,-30){\circle{20}}\put(180,-30){\makebox(0,0){12}}
\put(180,-46){\makebox(0,0){-}}
\put(180,46){\makebox(0,0){5}}
\put(180,20){\line(0,-1){40}}
\put(184,16){\makebox(0,0){4}}
\put(185,-16){\makebox(0,0){18}}
\put(240,0){\makebox(0,0){\normalsize $v=256$}}
\end{picture}$$
}

\medskip
Cases b), l), and m) suggest that for $K(3d,d)$ the group $G$
is elementary abelian of order $2^{3d-1}$, and this is indeed
easy to prove.\footnote{
Denote multiplication in $G$ by $*$, and let juxtaposition
denote disjoint union in the underlying $3d$-set $Z$. We show $A*B=B*A$.
Let $A = ER$, $B = ES$, where $E = A \cap B$.
Let $Z = EFGRST$, with $|E|=|F|=|G| = e$ and $|R|=|S|=|T|=d-e$.
Then $A*B = ER*ES = FS*FR = GR*GS = ES*ER = B*A$,
where $ER*ES = FS*FR$ since $ER*ES = FS*GT*GT*FR = FS*FR$,
and similarly for the other two equalities.
}
Thus, for each $d$ this approach yields only a
single graph $\Gamma$ that is locally $K(3d,d)$.

\end{document}